\documentclass[oneside,11pt]{amsart}
\usepackage{graphicx} 
\usepackage[margin=1in]{geometry}
\usepackage{amsmath, amssymb, amsfonts, amstext, amsthm, amscd}
\usepackage[mathscr]{euscript}
\usepackage{enumerate}
\usepackage{tikz,color,soul}
\usepackage[english]{babel}
\usepackage{chngcntr}
\newtheorem{theorem}{Theorem}[section]

\newtheorem{proposition}[theorem]{Proposition}

\theoremstyle{definition}
\newtheorem{definition}[theorem]{Definition}

\numberwithin{equation}{section}

\title{ON THE SOMBOR INDEX OF THE TOTAL GRAPH AND THE UNIT GRAPH OF COMMUTATIVE RINGS }
\author{ABHISHEK VAIBHAV PATHAK}
\address{St Joseph's University,36, Langford Rd, Langford Gardens, Bengaluru, Karnataka 560027}
\email{avpathak20@gmail.com}
\author{ANUKUL SACHAN}
\address{St Joseph's University,36, Langford Rd, Langford Gardens, Bengaluru, Karnataka 560027}
\email{anukulsachan5454@gmail.com}
\author{Raisa DSouza}
\address{St Joseph's University,36, Langford Rd, Langford Gardens, Bengaluru, Karnataka 560027}
\email{raisadsouza@sju.edu.in}
\date{}
\keywords{Sombor Index, Total Graph, Unit Graph, Local Rings, Topological Indices}
\subjclass[2020]{Primary 05C50}

\begin{document}
\maketitle
\begin{abstract}
In this paper, we investigate the Sombor index of the total graph and unit graph of $\mathbb{Z}_n$ which is denoted by $T_{\Gamma}(\mathbb{Z}_n)$ and $G(\mathbb{Z}_n)$ respectively for $n \in \{2k, p^{\alpha}, pq, p^2q\}$ where $p$ and $q$  are distinct odd prime numbers such that $p < q$. Moreover, we compute the Sombor index of any finite local ring.
\end{abstract}
\section{Introduction}
Algebra and Graph Theory are disciplines of mathematics that focus on building and investigating certain structures. In recent years, mathematicians have constructed many graphs from algebraic structures, considering different properties inherited by these structures. The properties of graphs have been studied to understand the algebraic structures better.

An edge connecting vertices $u$ and $v$ is denoted by $uv$ in a graph $G.$
In the mathematical and chemical literature,  many vertex–degree–based or distance based graph invariants (usually referred to as “topological indices”) have been introduced and
extensively studied \cite{two}. The general formula of a vertex-degree-based topological index is,
\begin{equation*}
    TI(G)=\sum_{uv \in E(G)} H(d_u,d_v)
\end{equation*}
where $H(x,y)$ is some real valued function with the property $H(x,y)=H(y,x)$ and $d_x$ is the degree of vertex $x.$ Similarly the general formula for a distance based topological index is,
\begin{equation*}
    TI(G)=\sum_{\{u,v\}\subseteq V(G)} H(d(u,v))
\end{equation*}
where $d(u,v)$ is the distance between the vertices $u$ and $v.$
A topological index of $G$ is a real number associated with $G.$  It is independent of a graph's labeling (or pictorial representation) of a graph. 
In 2021 mathematical chemist Ivan Gutman introduced the Sombor Index \cite{twelve} of a graph as defined below 
$$SO(G)=\sum_{uv\in E(G)} \sqrt{d_u^2 + d_v^2.}$$
Further in the same year Saeid Alikhani and Nina Ghanbari \cite{five} gave the formulae for the Sombor Index of cycle graphs, bipartite graphs, and wheel graphs. They also presented theorems concerning the aforementioned index, which are utilized to generalize the Sombor Index for various graphs, including those with ring structures. In 2022, Arif, Alper, and Necla \cite{arif@2022} studied the Sombor index of zero-divisor graphs of $\mathbb{Z}_n$. 

In this paper, we study the Sombor index of the total graph and unit graph of some commutative rings. In Section 2, we recall some basic definitions, notations, and results which will be used throughout the paper. Also, we calculate the Sombor index of the total graph and the unit graph of $\mathbb{Z}_n$ in Section 3 and Section 4 respectively. In Section 5, we explore the Sombor index the total graph, and the unit graph of any finite local ring. We establish a relationship between the complete graph, $k$ -  regular graph, and its complement in Section 6.

\section{PRELIMINARIES}
In this section, we recall some basic definitions, notations, and results which will be used throughout the paper. We will assume all rings are commutative with unity. Let $Z(R)$ denote the set of zero-divisors and $U(R)$ denote the set of units of the ring $R.$

\begin{definition}
    The total graph \cite{three} of the ring $R,$ denoted $T_\Gamma(R),$ is the graph obtained by setting all the elements of $R$ to be the vertices and distinct vertices $x$ and $y$ are adjacent if and only if $x + y \in Z(R).$
\end{definition}
The degree of each vertex in a total graph can easily be determined. It was given by T.Asir and Thirugnanam Tamizh Chelvam \cite{AST}, among other results. We state the results below for completeness. 
\begin{theorem}\cite{AST}
\label{P4}
    Let $R$ be a ring. The following statements hold for the total graph of $R$
            \begin{enumerate}
                \item If $2 \notin U(R)$, $d_x = n-\left|U(R)\right|-1$, for every $x \in R$
                \item If $2 \in U(R)$, $d_x = n-\left|U(R)\right|-1$, for every $x \in Z(R)$ and $d_x = n-\left|U(R)\right|,$ for every $x \notin Z(R).$         
            \end{enumerate}          
\end{theorem}
\begin{theorem}\cite[Theorem 2.5]{AST}
\label{P2}
    The following statements hold for the total graph of $\mathbb{Z}_n:$
            \begin{enumerate}
                \item If $n$ is even then $d_v = n-\phi(n)-1,$ for every $v \in \mathbb{Z}_n.$
                \item If $n$ is odd then $d_v = n-\phi(n)-1,$ for every $v \in Z(\mathbb{Z}_n)$ and $d_v = n-\phi(n),$ for every $v \notin Z(\mathbb{Z}_n),$ where, $\phi$ is the Euler function.         
            \end{enumerate}          
\end{theorem}

\begin{definition}
The unit graph \cite{one} of $R,$ denoted $G(R),$ is the graph obtained by setting all the elements of $R$ to be the vertices and  distinct vertices $x$ and $y$ are adjacent if and only if $x + y \in U(R).$
\end{definition}
\newpage
Analogous to Theorem \ref{P4}, in \cite{article1} gave the degrees of vertices in $G(R).$
\begin{theorem}
\label{P7}
    Let $R$ be a ring. The following statements hold for the unit graph of $R$
            \begin{enumerate}
                \item If $2 \notin U(R)$, $d_x = \left|U(R)\right|$, for every $x \in R$
                \item If $2 \in U(R)$, $d_x = \left|U(R)\right|-1$, for every $x \in Z(R)$ and $d_x = \left|U(R)\right|,$ for every $x \notin Z(R).$         
            \end{enumerate}          
\end{theorem}
\begin{theorem}\cite[Proposition 2.4]{UG@2010}
\label{P1}
     Then the following statements hold for the unit graph of $\mathbb{Z}_n:$ 
         \begin{enumerate}
             \item If $2 \notin U(\mathbb{Z}_n),$ then $d_x = \phi(n)$ for every $x \in \mathbb{Z}_n.$
             \item If $2 \in U(\mathbb{Z}_n),$ then $d_x = \phi(n)-1$ for every $x \in U(\mathbb{Z}_n)$ and $d_x = \phi(n)$ for every $x \in \mathbb{Z}_n\backslash U(\mathbb{Z}_n).$
         \end{enumerate}
\end{theorem}

\begin{proposition}\cite{Unitgrah}\label{P6}
     Let $R$ be a finite ring, then $Z(T_\Gamma(R))$ is a complete graph if and only if $R$ is a local ring.
\end{proposition}

\section{SOMBOR INDEX OF THE TOTAL GRAPH OF $\mathbb{Z}_n$}
Recently, the total graph of the ring $\mathbb{Z}_n$
  has gained popularity in spectral and chemical graph theory. Many researchers have focused on this area, contributing to its advancement and exploration. Sheela Suthar and Om Prakash have examined the energy and Wiener index of the total graph $T_{\Gamma}(\mathbb{Z}_n$) \cite{SUTHAR2017485}. In this section, we analyze the Sombor index the total graph of $\mathbb{Z}_n$.
\begin{theorem}
     Let $p$ be a prime number.The following holds for the total graph of $\mathbb{Z}_n:$
        \begin{enumerate}
            \item If $n$ is even, then $SO(T_\Gamma(\mathbb{Z}_n)) = \frac{n(n-\phi(n)-1)^2}{\sqrt{2}}.$
            \item If $n=p^\alpha,$ where $\alpha \geq 1$ and $p > 2,$ then 
            \begin{align*}
                  SO(T_\Gamma(\mathbb{Z}_n)) = \frac{\phi(n)(n-\phi(n))^2}{\sqrt{2}} + \frac{(n-\phi(n)-1)^2(n-\phi(n))}{\sqrt{2}}.
            \end{align*} 
           \end{enumerate}
\end{theorem}
\begin{proof}
  \begin{enumerate}
      \item If $n$ is even, $T_\Gamma(\mathbb{Z}_n)$ is a $n-\phi(n)-1$ regular graph with $n$ vertices and the result follows. 
      \item For $n=p^\alpha,$ by Theorem \ref{P2}, if $x \in Z(\mathbb{Z}_n)$ then,  $d_x=n-\phi(n)-1$  and if $x \notin Z(\mathbb{Z}_n)$ then, $d_x=n-\phi(n).$ Now, we know that for $n=p^\alpha$ the set of non-units of $\mathbb{Z}_n$ forms a subgroup. Hence, any two non-units of $\mathbb{Z}_n$ are adjacent. Therefore, number of edges between zero-divisors is, $$\frac{(n-\phi(n))(n-\phi(n)-1)}{2}.$$
            By handshake lemma, the total number of edges in $T_\Gamma(\mathbb{Z}_n)$ is given by, $$\frac{\phi(n)(n-\phi(n)) + (n-\phi(n))(n-\phi(n)-1)}{2}.$$ Thus, the number of edges among units will be, 
                $\displaystyle{\frac{\phi(n)(n-\phi(n))}{2}}.$
                The result follows.   
  \end{enumerate}  
\end{proof}
\begin{theorem}
Let $p$ and $q$ be two odd prime numbers with $p < q.$ Then, the Sombor index of the graph $T_\Gamma(\mathbb{Z}_{pq})$ is given by,
    \begin{multline*}
    SO(T_\Gamma(\mathbb{Z}_{pq})) = \sqrt{2}\alpha(pq-\phi(pq)-1) + \\
\beta\sqrt{(pq-\phi(pq)-1)^2+(pq-\phi(pq))^2} +\\ \sqrt{2}(\left|E\right|-\alpha-\beta)(pq-\phi(pq))
    \end{multline*}
    where, \begin{align*}   
    \alpha = \frac{p(p-1)+q(q-1)}{2},\quad \beta = 2(p-1)(q-1),
    \quad \left|E\right| = \frac{(pq-1)(p+q-1)}{2}.
    \end{align*} 
\end{theorem}
\begin{proof}
    Let $Z(\mathbb{Z}_{pq}) = S_p \cup S_q$ where,  $S_p = \{0,p,2p,\cdots,(q-1)p\}$ and $S_q = \{0,q,2q,\cdots,(p-1)q\}.$ We observe that $S_p$ and $S_q$ are closed under the operation $\bigoplus_{pq}.$ The fact that $S_p$ and $S_q$ are closed and $S_p\cap S_q = {0}$, will imply that $a \nsim b \,$ for all $ \,  a \in S_p \, \text{and} \, b \in S_q.$ We thus get complete graphs $K_q$ and $K_p$ for the sets $S_p$ and $S_q$ respectively. Thus the number of edges between zero divisors and zero divisors are $$\alpha = \frac{q(q-1)}{2} + \frac{p(p-1)}{2}.$$ Since $p$ and $q$ are odd primes, by Theorem \ref{P2} $d_v = pq - \phi(pq) - 1, \,$ for all $\, v \in Z(\mathbb{Z}_{pq}).$ Thus,
    $$\begin{aligned}
        d_v & =  pq - \phi(pq) - 1 \\
        & = pq - (p-1)(q-1) - 1 \\
        & = pq - pq + p + q - 1 - 1 \\
        & = p + q - 2.
    \end{aligned}$$
    We notice that the degree of a vertex within $S_p$ is, $d^{S_p}_v = q-1 \,$ for all $ \, v \in S_p$ and within $S_q$ is, $d^{S_q}_v = p-1 \, $ for all$\, v \in S_q.$ So the remaining adjacencies of $v \in S_p$ or $\, S_q \,$ in $T_\Gamma(\mathbb{Z}_{pq})$ is with the units. The excess degree of $$v \in S_p \, \text{ is } \, p + q - 2 - q + 1 = p - 1$$ and of $$v \in S_q \,\text{ is } \, p + q - 2 - p + 1 = q - 1.$$ Next we notice that $0 \nsim u \, $ for all$ \, u \in U(\mathbb{Z}_{pq}).$ Thus the number of edges between zero divisors and units are $$ \begin{aligned} 
    \beta & = (q-1)(p-1) + (p-1)(q-1) \\ 
    & = 2(p-1)(q-1).
    \end{aligned}$$
    Using Theorem \ref{P2} and the handshake lemma we find the total number of edges in $T_\Gamma(\mathbb{Z}_{pq})$ as, 
    $$\begin{aligned}
        \left|E\right|  = \frac{1}{2}(pq-1)(p+q-1).
    \end{aligned}$$ 
    Since we know the total number of edges in $T_\Gamma(\mathbb{Z}_{pq}) \, \,(\left|E\right|)$, the number of edges between zero divisors and zero divisors $(\alpha)$ and the number of edges between zero divisors and units $(\beta),$ we can easily get the number of edges between units and units, which will be given by $\left|E\right| - \alpha - \beta.$
    
    From Theorem \ref{P2} we know $d_v = pq-\phi(pq)-1 \, \text{for all} \, v \, \in Z(\mathbb{Z}_{pq})$ and  $d_u = pq-\phi(pq) \,$ for all $\, u \, \in U(\mathbb{Z}_{pq}).$ The result follows. 
\end{proof}

\begin{theorem}
Let $p$ and $q$ be two odd prime numbers with $p < q.$ Then, the Sombor index of the graph $T_\Gamma(\mathbb{Z}_{p^2q})$ is 
    \begin{multline*}
    SO(T_\Gamma(\mathbb{Z}_{p^2q})) = \sqrt{2}\alpha(p^2q-\phi(p^2q)-1) + \\
\beta\sqrt{(p^2q-\phi(p^2q)-1)^2+(p^2q-\phi(p^2q))^2} +\\ \sqrt{2}(\left|E\right|-\alpha-\beta)(p^2q-\phi(p^2q))
    \end{multline*}
    where, \begin{align*}
         \alpha& =  \frac{p(q-1)(p(q-1)-1)}{2}+ \frac{p(p-1)(p(p-1)-1)}{2}
         \\ & \hspace{4.6cm} +\frac{p(p-1)}{2} +p^2(q-1)+p^2(p-1), \\  \beta& =  2p^2
    (p-1)(q-1), \, \\ 
    \left|E\right|& =  \frac{p(p+q-1)(p^2q-1)}{2}.
    \end{align*}  
\end{theorem}
\begin{proof}
    Let $Z(\mathbb{Z}_{pq}) = S_p \cup S_q \cup S_{pq}$ where, $S_p = \{0,p,2p,\cdots,(pq-1)p\}$, $S_q = \{0,q,2q,\cdots,(p^2-1)q\} $ and $ S_{pq}= \{0,pq,2pq,\cdots,(p-1)pq\}.$ Now consider $S^*_p = S_p\backslash S_{pq}$ and $S^*_q = S_p\backslash S_{pq}.$ We observe that $S_{pq},$ $S^*_p$ and $S^*_q$ are closed under $\bigoplus_{p^2q}.$ The fact that $S^*_p$ and $S^*_q$ are closed and $S^*_p\cap S^*_q = \phi$, will imply that $a \nsim b \, $ for all$ \,  a \in S^*_p \, \text{and} \, b \in S^*_q.$ We thus get complete graphs $K_p,$ $K_{p(q-1)}$ and $K_{p(p-1)}$ for the sets $S_{pq},\, S^*_p, \, S^*_q \,$ respectively. We also observe that $a \sim b \, $ for all$ \,  a \in S^*_p \, \text{and} \, b \in S_{pq}$ and $a \sim b \, $ for all$ \,  a \in S^*_q \, \text{and} \, b \in S_{pq}.$ We thus get complete bipartite graphs $K_{p(q-1),p} \,$ between the vertices of $S^*_p$ and $S_{pq}$ and $K_{p(p-1),p} \,$ between the vertices of $S^*_q$ and $S_{pq}$. Therefore the number of edges between zero-divisors and zero-divisors are, $$\alpha = \frac{p(q-1)(p(q-1)-1)}{2} + \frac{p(p-1)(p(p-1)-1)}{2} + \frac{p(p-1)}{2} + p^2(q-1) + p^2(p-1).$$
    Since $p$ and $q$ are odd primes, by Theorem \ref{P2} $d_v = p^2q - \phi(p^2q) - 1, \,$ for all $\, v \in Z(\mathbb{Z}_{p^2q}).$ Thus,
    $$\begin{aligned}
        d_v & =  p(p + q - 1) - 1.
    \end{aligned}$$
     We notice that $d^{S^*_p}_v = p(q-1) -1+p \,$ for all $ \, v \in S_p$ and  $d^{S^*_q}_v = p(p-1) -1 +p \, $ for all$\, v \in S_q.$ So the remaining adjacencies of $v \in S^*_p$ or $\, S^*_q \,$ in $T_\Gamma(\mathbb{Z}_{p^2q})$ is with the units. The excess degree of $$v \in S^*_p \, \text{ is } \,p(p + q - 1) - 1 - (p(q-1) -1+p) = p(p - 1)$$ and of $$v \in S^*_q \,\text{ is } \,p(p + q - 1) - 1 -(p(p-1) -1 +p) = p(q - 1).$$ Also the excess degree of  $v \in S_{pq} \, is \,0$ as degree of any element in $S_{pq}$ is $p^2q - \phi(p^2q) - 1 .$ Thus the number of edges between zero-divisors and units are, $$\begin{aligned}
         \beta & = p(p-1)p(q-1) + p(q-1)p(p-1) \\ & = 2p^2(p-1)(q-1).
         \end{aligned}$$
         Using Theorem \ref{P2} and the handshake lemma we find the total number of edges in $T_\Gamma(\mathbb{Z}_{p^2q})$ as follows, 
    $$\begin{aligned}
        \left|E\right| & = \frac{1}{2}[\phi(p^2q)(p^2q-\phi(p^2q)) + (p^2q-\phi(p^2q))(p^2q-\phi(p^2q)-1]) \\
        & = \frac{1}{2}p(p+q-1)(p^2q-1).
    \end{aligned}$$ 
    Since we know the total number of edges in $T_\Gamma(\mathbb{Z}_{pq}) \, \,(\left|E\right|)$, the number of edges between zero divisors and zero divisors $(\alpha)$ and the number of edges between zero divisors and units $(\beta),$ we can easily get the number of edges between units and units, which will be given by $\left|E\right| - \alpha - \beta.$
    
    From Theorem \ref{P2} we know $d_v = pq-\phi(pq)-1 \, \text{for all} \, v \, \in Z(\mathbb{Z}_{pq})$ and  $d_u = pq-\phi(pq) \,$ for all $\, u \, \in U(\mathbb{Z}_{pq}).$ The result follows.

\end{proof}
\section{SOMBOR INDEX OF UNIT GRAPH OF $\mathbb{Z}_n$}
An intriguing observation lies in the relationship between unit graphs and total graphs: they are complements of each other. In the preceding section, we delved into computing the Sombor index of total graphs. Now, in this section, we focus on exploring the Sombor index of unit graphs.
\begin{theorem}
Let $p$ be a prime number. The following holds for unit graph of $\mathbb{Z}_n:$
        \begin{enumerate}
            \item If $n$ is even, then $SO(G(\mathbb{Z}_n)) = \frac{n(\phi(n))^2}{\sqrt{2}}.$
            \item If $n=p^\alpha,$ where $\alpha \geq 1$ and $p > 2$ then 
\begin{multline*}
            SO(G(\mathbb{Z}_n)) = \phi(n)(n-\phi(n))\sqrt{(\phi(n))^2 + (\phi(n)-1)^2} + \\ \frac{[\phi(n)(\phi(n)-1) - (n-\phi(n))](\phi(n)-1)}{\sqrt{2}.}
            \end{multline*}
           \end{enumerate}
\end{theorem}
\begin{proof}
\begin{enumerate}
            \item For $n$ even, as $2 \notin U(\mathbb{Z}_n),$ the degree of each vertex is $\phi(n)$ by Theorem \ref{P1}. Thus, $G(\mathbb{Z}_n)$ is a $\phi(n)$-regular graph with $\frac{n\phi(n)}{2}$ edges. The result follows.
           \item For $n=p^\alpha, 2\in U(\mathbb{Z}_n).$ By Theorem \ref{P1},  if $x \in U(\mathbb{Z}_n)$, then $d_x=\phi(n)-1,$  and if $x \notin U(\mathbb{Z}_n)$ then,  $d_x=\phi(n).$ Now, we know that for $n=p^\alpha$ the set of non-units of $\mathbb{Z}_n$ forms a group under modulo addition. Hence, any two non-units of $\mathbb{Z}_n$ are not adjacent. But, we know that degree of each non-unit is $\phi(n)$, thus number of edges between units and non-units is $(n-\phi(n))\phi(n).$
            Now by Handshake Lemma, total number of edges in $G(\mathbb{Z}_n)$ is given by $$\frac{\phi(n)(\phi(n)-1) + (n-\phi(n))\phi(n)}{2}.$$ Thus, the number of edges among units will be the difference between the total number of edges and the number of edges among units and non-units,i.e. 
                $$\frac{\phi(n)(\phi(n)-1) - (n-\phi(n))\phi(n)}{2} = \frac{\phi(n)(\phi(n)-1-n-\phi(n))}{2}.$$ 
                The result follows.
\end{enumerate}
\end{proof}
\begin{theorem}
    Let $p$ and $q$ be two odd prime numbers with $p < q.$ Then, Sombor index of graph $G(\mathbb{Z}_{pq})$ is,
    $$SO(G(\mathbb{Z}_{pq})) = \sqrt{2}\alpha\phi(pq) +
    \beta\sqrt{(\phi(pq))^2+(\phi(pq)-1)^2} + \sqrt{2}(\left|E\right|-\alpha-\beta)(\phi(pq)-1)$$
    where,  
    $\alpha = (p-1)(q-1), \quad \beta = (p-1)(q-1)(p+q-3),
    \quad \left|E\right| = \displaystyle{\frac{(pq-1)\phi(pq)}{2}}.$
    
\end{theorem}
\begin{proof}
    Let $Z(\mathbb{Z}_{pq}) = S_p \cup S_q$ where $S_p = \{0,p,2p,\cdots,(q-1)p\}$ and $S_q = \{0,q,2q,\cdots,(p-1)q\}.$ We can observe that the sets $S_p$ and $S_q$ forms a closed set under the operation $\bigoplus_{pq}.$ Since the unit graph is the complement of the total graph and in the total graph, $S_p$ and $S_q$ induce complete graphs, thus there will be no edges inside $S_p$ and $S_q$. But we will get a complete bipartite graph with the partite sets $S_p$ and $S_q$. 
    
    Also observe that $0 \nsim v \, $ for all$ \, v \, \in Z(\mathbb{Z}_{pq})$. Thus the sets $S_p$ and $S_q$ gives us $K_{(p-1),(q-1)}.$ So that the number of edges between zero divisors and zero divisors is, $$\alpha = (p-1)(q-1).$$
    Since $p$ and $q$ are odd primes, then by Theorem \ref{P1} $d_v = \phi(pq)  \, $ for all$ \, v \in Z(\mathbb{Z}_{pq}).$ Thus,
    $$\begin{aligned}
        d_v & = \phi(pq) = (p-1)(q-1)
    \end{aligned}$$
    Notice that for all $ \, v \in Z(\mathbb{Z}_{pq}) \, d_v=\phi(pq) = (p-1)(q-1).$ Also, for all $\, v \in S_p\backslash\{0\} \, d_v=(p-1).$ and for all $\, v \in S_q\backslash\{0\} \, d_v=(q-1).$ 
    Thus the remaining degrees of, $$v \in S_p\backslash\{0\} \, \text{ is } \, (p-1)(q-1)- (p-1)=(p-1)(q-2)$$ and of $$v \in S_q\backslash\{0\} \, \text{ is } \, (p-1)(q-1)- (q-1)=(p-2)(q-1).$$ Also $0 \sim u\,$ for all $\, u \, \in U(\mathbb{Z}_{pq})$, which implies $d_0=\phi(pq)=(p-1)(q-1).$ Therefore the number of edges between zero divisors and units are, $$ \begin{aligned} 
    \beta & = (q-1)(p-1)(q-2) + (p-1)(p-2)(q-1) + (p-1)(q-1) \\ 
    & = (p-1)(q-1)(p+q-3)
    \end{aligned}$$
    Using Theorem \ref{P1} and handshake lemma we can find the total number of edges in $G(\mathbb{Z}_{pq})$ as follows, 
    $$\begin{aligned}
        \left|E\right| & = \frac{1}{2}[\phi(pq)(pq-\phi(pq)) + (\phi(pq)-1)(\phi(pq)]) \\
        & = \frac{1}{2}((p-1)(q-1))(pq-1).
    \end{aligned}$$ 
    Now since we know the total number of edges in $G(\mathbb{Z}_{pq})$ $(\left|E\right|)$, number of edges between zero divisors and zero divisors $(\alpha)$ and the number of edges between zero divisors and units $(\beta)$, we can easily get the number of edges between units and units, which will be given by $\left|E\right| - \alpha - \beta.$
    
    From Theorem \ref{P1} we know $d_v = \phi(pq) \,$ for all $\, v \, \in Z(\mathbb{Z}_{pq})$ and $d_u = \phi(pq)-1 \, $ for all $ \, u \, \in U(\mathbb{Z}_{pq}).$ The result follows. 
\end{proof}

\begin{theorem}
Let $p$ and $q$ be two odd prime numbers with $p < q.$ Then, the Sombor index of the graph $G(\mathbb{Z}_{p^2q})$ is,
    $$SO(G(\mathbb{Z}_{p^2q})) = \sqrt{2}\alpha\phi(p^2q) + 
\beta\sqrt{(\phi(p^2q))^2+(\phi(p^2q)-1)^2} + \sqrt{2}(\left|E\right|-\alpha-\beta)(\phi(p^2q)-1)$$
    where, $
        \alpha =  p^2(p-1)(q-1),\quad \beta = p^2(p-1)(q-1)(p+q-3), \quad
    \left|E\right| = \displaystyle{\frac{p^2(p-1)(q-1)(p^2q-1)}{2}.}$
\end{theorem}
\begin{proof}
Let $Z(\mathbb{Z}_{p^2q}) = S_p \cup S_q \cup S_{pq}$ where, $S_p = \{0,p,2p,\cdots,(pq-1)p\}$,  $S_q = \{0,q,2q,\cdots,(p^2-1)q\} $ and $ S_{pq}= \{0,pq,2pq,\cdots,(p-1)pq\}.$ Now consider  $S^*_p = S_p\backslash S_{pq}$ and $S^*_q = S_p\backslash S_{pq}.$ Since the unit graph is the complement of the total graph and in the total graph, $S^*_p$, $S^*_q$ and $S_{pq}$ induce complete graphs, thus here there will be no edges inside $S^*_p$, $S^*_q$ and $S_{pq}.$ Also, in the total graph we get complete bipartite graphs induced by $S^*_p$ \& $S_{pq}$ and $S^*_q$ \& $S_{pq}.$ Thus here there will be no edges between $S^*_p$ \& $S_{pq}$ and  $S^*_q$ \& $S_{pq}.$ So the only edges between zero-divisors and zero-divisors will be between $S^*_p$ and $S^*_q$, which will give us a complete bipartite graph. Also observe that $0 \nsim v \, $ for all$ \, v \, \in Z(\mathbb{Z}_{p^2q})$.Thus the number of edges between zero divisors and zero divisors are, $$\alpha = p^2(p-1)(q-1).$$
    Since $p$ and $q$ are odd primes, then by Theorem \ref{P1} $d_v = \phi(p^2q) \, $ for all$ \, v \in Z(\mathbb{Z}_{p^2q}).$ Thus,
    $$\begin{aligned}
        d_v & = \phi(p^2q) = p(p-1)(q-1).
    \end{aligned}$$
    Notice that for all $ \, v \in Z(\mathbb{Z}_{p^2q}) \, d_v=\phi(p^2q) = p(p-1)(q-1).$ Also, for all $\, v \in S_p^* \, d_v=p(q-1).$ and for all $\, v \in S_q^* \, d_v=p(p-1).$ 
    Thus the remaining degrees of $$v \in S_p^* \, \text{is} \, p(p-1)(q-1)- p(q-1)= p(p-2)(q-1)$$ and of $$v \in S_q^* \, \text{is} \, p(p-1)(q-1)- p(p-1)=p(p-1)(q-2).$$ Also each element in $S_{pq}$ is adjacent to every unit. Therefore the number of edges between zero divisors and units are, $$ \begin{aligned} 
    \beta & = p^2(q-1)(p-1)(p-2) + p^2(p-1)(q-2)(q-1) + p^2(p-1)(q-1) \\ 
    & = p^2(p-1)(q-1)(p+q-3)
    \end{aligned}$$
    Using Theorem \ref{P1} and the handshake lemma we can find the total number of edges in $G(\mathbb{Z}_{p^2q})$ as follows,
    $$\begin{aligned}
        \left|E\right| & = \frac{1}{2}[\phi(p^2q)(p^2q-\phi(p^2q)) + (\phi(p^2q)-1)(\phi(p^2q)]) \\
        & = \frac{1}{2}(p^2(p-1)(q-1))(p^2q-1).
    \end{aligned}$$ 
    Now since we know the total number of edges in $G(\mathbb{Z}_{p^2q})$ $(\left|E\right|)$, the number of edges between zero divisors and zero divisors $(\alpha)$ and the number of edges between zero divisors and units $(\beta)$, we can easily get the number of edges between units and units, which will be given by $\left|E\right| - \alpha - \beta.$
    
    From Theorem \ref{P1} we know $d_v = \phi(p^2q) \,$ for all $\, v \, \in Z(\mathbb{Z}_{p^2q})$ and $d_u = \phi(p^2q)-1 \, $ for all $ \, u \, \in U(\mathbb{Z}_{p^2q}).$ The result follows.
\end{proof}
\section{SOMBOR INDEX OF LOCAL RINGS}
It is easy to calculate the Sombor index of the total graph and the unit graph of the rings $\mathbb{Z}_{p^\alpha},$ where $p$ is a prime number and $\alpha$ is a natural number. In this section, we generalize the computation of the Sombor index of the total graph and the unit graph of any finite local ring.
\begin{theorem} \label{LR}
    Let $R$ be a local ring
    \begin{enumerate}
                \item If $2 \notin U(R)$ then, $$SO(T_\Gamma(R) = \frac{n(n-\left|U(R)\right|-1)^2}{\sqrt{2}}$$
                \item If $2 \in U(R)$ then, $$SO(T_\Gamma(R)) = \frac{(n-\left|U(R)\right|)(n-\left|U(R)\right|-1)^2}{\sqrt{2}}+\frac{(\left|U(R)\right|)(n-\left|U(R)\right|)^2}{\sqrt{2}}.$$         
            \end{enumerate}
\end{theorem}
\begin{proof}
    \begin{enumerate}
        \item If $2 \notin U(R)$, then by Theorem \ref{P4} we know $d_x=n-\left|U(R)\right|-1$, for every $x \in R.$ Thus by the handshake lemma we get the total number of edges as $$\left| E\right|=\frac{n(n-\left|U(R)\right|-1)}{2}$$ 
        and the result follows.
        \item If $2 \in U(R)$, then by Theorem \ref{P6} $Z(T_\Gamma(R))\cong K_{\left|Z(R)\right|}$, thus the number of edges between zero divisors and zero divisors are given by $$\left|E_1\right| = \frac{\left|Z(R)\right|(\left|Z(R)\right|-1)}{2}=\frac{(n-\left|U(R)\right|)(n-\left|U(R)\right|-1)}{2}$$
        By Theorem \ref{P4} when $2 \in U(R)$ we know that $d_x=n-\left|U(R)\right|-1$, for every \\$x \in Z(R).$ Also $\left|Z(R)\right| = n -\left|U(R)\right|$ and  $Z(T_\Gamma(R))\cong K_{\left|Z(R)\right|}$ which implies \\ $d_x=n-\left|U(R)\right|-1$, for every $x \in Z(R)$ i.e. all the degrees of zero divisors are exhausted. Thus there will be no edges between zero divisors and units.\\ From Theorem \ref{P4} when $2 \in U(R)$ we know $d_x=n-\left|U(R)\right|$, for every $x \notin Z(R).$ Thus the number of edges between units and units is given by, $$\left|E_2\right| = \frac{(n-\left|U(R)\right|)(n-\left|U(R)\right|-1)}{2}.$$ Again the result follows.      
    \end{enumerate}
\end{proof}
\begin{theorem}
    Let $R$ be a local ring
    \begin{enumerate}
                \item If $2 \notin U(R)$ then, $$SO(G(R) = \frac{n(\left|U(R)\right|)^2}{\sqrt{2}}$$
                \item If $2 \in U(R)$ then, $$SO(G(R))= {\left|U(R)\right|\left(n-\left|U(R)\right|\right)\sqrt{\left|U(R)\right|^2+(n-\left|U(R)\right|)^2}}$$       
            \end{enumerate}
\end{theorem}
\begin{proof}
    The proof follows similar to Theorem \ref{LR}.
\end{proof}
\section{SOMBOR INDEX OF $k$- REGULAR GRAPH AND ITS COMPLEMENT}
In this section, we give the relationship between the Sombor index of a complete graph, $k$-regular graph, and its complement.
\begin{theorem}
    Let $G$  be a $k - $  regular graph of order $n$ and $\Bar{G}$ be its complement.  Then,
    \begin{equation*}
        SO(K_n) = \left(\sqrt{SO(G)} + \sqrt{SO(\Bar{G})}\right)^2
    \end{equation*}
\end{theorem}
\begin{proof}
    As the graph $G$ is a $k -$ regular graph, it is clear that the graph $\Bar{G}$ is a $(n-k-1)-$ regular graph. Therefore, 
    \begin{equation*}
        SO(G) = \frac{nk^2}{\sqrt{2}}
\hspace{1cm} and \hspace{1cm}
        SO(\Bar{G}) = \frac{n(n-k-1)^2}{\sqrt{2}}.
    \end{equation*}
Now, we know that, 
$$\begin{aligned}
    SO(K_n) & = \frac{n(n-1)^2}{\sqrt{2}} \\
    & = \frac{n(k + (n-k-1))^2}{\sqrt{2}} \\
    & = SO(G) + \sqrt{2}n\sqrt{\left(\frac{\sqrt{2}SO(G)}{n}\right)}\sqrt{\left(\frac{\sqrt{2}SO(\Bar{G})}{n}\right)} + SO(\Bar{G}) \\
     & = \left(\sqrt{SO(G)} + \sqrt{SO(\Bar{G})}\right)^2
\end{aligned}$$
\end{proof}
\bibliography{references.bib}
\bibliographystyle{ieeetr}
\end{document}